\documentclass[a4paper]{amsproc}
\usepackage{amssymb}
\usepackage{amscd} %% Package for commutative diagrams
\lccode`\-=`\-
\newcommand{\grad}{\mathop{\rm grad}\nolimits}
\renewcommand{\div}{\mathop{\rm div}\nolimits}

\theoremstyle{plain}
 \newtheorem{thm}{Theorem}[section]

\theoremstyle{definition}

\numberwithin{equation}{section}

\renewcommand{\leq}{\leqslant}
\renewcommand{\geq}{\geqslant}

\title[Domain decomposition schemes for the Stokes equation]
{Domain decomposition schemes for the Stokes equation}


\author[Vabishchevich]{\bfseries Petr N. Vabishchevich}

\address{
Keldysh Institute of Applied Mathematics \\ % \hfill (Received 00 00 2010)\\
4 Miusskaya Square  \\ %\hfill (Revised  00 00 2010)\\
125047 Moscow\\ 
Russia}
\email{vabishchevich@gmail.com}

\dedicatory{Dedicated to Academician Anton Bilimovic}
\begin{document}

%{\begin{flushleft}\baselineskip9pt\scriptsize
%PUBLICATIONS DE L'INSTITUT MATH\'EMATIQUE\newline
%Nouvelle s\'erie, tome 87(101) (2010), od--do \hfill DOI:
%\end{flushleft}}
\vspace{18mm}
\setcounter{page}{1}
\thispagestyle{empty}

\begin{abstract}
Numerical  algorithms for solving  problems  of  mathematical
physics  on  modern  parallel  computers  employ various  domain
decomposition techniques.  Domain decomposition schemes are developed here
to solve numerically initial/boundary value  problems
for   the  Stokes  system of equations    in   the    primitive   variables
pressure-velocity. Unconditionally stable schemes of domain decomposition are based 
on the partition of unit for a computational domain and the corresponding 
Hilbert spaces of grid functions.
\end{abstract}

\maketitle

\section{Introduction}  

In computational fluid dynamics 
\cite{anderson1995computational,temam2001navier} 
there are employed numerical algorithms based on using the primitive
variables pressure-velocity. The main difficulties in this approach are
connected with the calculation of the pressure. In studying transient problems
the corresponding elliptic Neumann problem for the pressure is derived as
the result of employment of one or another scheme of splitting with
respect to physical processes \cite{glowinski2003finite,guermond2006overview}.

Domain decomposition methods are used for the numerical 
solution of boundary value problems 
for partial differential equations on parallel computers.  
They are in most common use for stationary 
problems  \cite{quarteroni1999domain,toselli2005domain}.
Computational algorithms with and without overlapping of 
subdomains are employed in this case 
in synchronous (sequential) and asynchronous (parallel) algorithms.  

For transient problems it seems to be more suitable to utilize iteration-free
variants of domain decomposition techniques 
\cite{mathew2008domain,samarskii2002difference} 
which are best suited to peculiarities of a problem (evolution in time). In these
regionally-additive schemes a transition to a new time level is performed via
solving problems in particular subdomains.

The regionally-additive schemes for the
Navier-Stokes equations in the primitive variables are discussed
in \cite{Churbanov2004domain}. In simulation of
incompressible flows an elliptic problem for the pressure can be changed to
separate elliptic problems for the pressure in particular subdomains. Therefore,
it is possible to construct iteration-free regionally-additive schemes for the
Navier-Stokes equations.
In this paper we propose a general approach to construct
domain decomposition schemes for time-dependent systems of equations. 
Using the partition of unit for a computational domain and 
the corresponding Hilbert spaces of grid functions we perform 
a transition to finding the individual components of the solution 
in the subdomains. The unsteady Stokes equations for an incompressible fluid is considered 
as a typical problem.

\section{Stokes equations}

Assume that the linear approximation is valid to describe a flow of an incompressible viscous fluid.
In a region $\Omega$ with solid boundaries  we can write equations of motion and continuity 
in the primitive variables \textit {pressure, velocity} as follows 
\begin{equation}\label{2.1}
  \frac{\partial {\bf u}}{\partial t}
 + {\rm grad} \, p
 - \nu ~\Delta \, {\bf u}
 = {\bf f} ( {\bf x},t),  
\end{equation}
\begin{equation}\label{2.2}
   {\rm div} \, {\bf u} = 0,
   \quad  {\bf x} \in \Omega, \quad 0 < t \leq T.  
\end{equation}
Here $\bf u$ is the velocity, $p$ is the pressure, $\nu$ is the kinematic viscosity 
and $\Delta = {\rm div}~{\rm grad}$ is the Laplace operator. 
Equations (\ref{2.1}), (\ref{2.2}) are supplemented with the following condition 
for the single-valued  evaluation of the pressure 
\begin{equation}\label{2.3}
   \int \limits_{\Omega} p( {\bf x},t) d {\bf x} = 0,
    \quad  0 < t \leq T.
\end{equation}
No-slip, no-permeability conditions are specified on solid boundaries 
\begin{equation}\label{2.4}
    {\bf u} ( {\bf x},t) = 0,
    \quad   {\bf x} \in \partial \Omega, \quad  0 < t \leq T.
\end{equation}
Some initial condition is also given
\begin{equation}\label{2.5}
    {\bf u} ({\bf x},0) =  {\bf v}({\bf x}),
    \quad   {\bf x} \in \Omega.
\end{equation}

Let us rewrite problem (\ref{2.1})--(\ref{2.5}) in an operator formulation. 
On the set of functions satisfying (\ref{2.3}), \ref{2.4}), we have the Cauchy problem 
\begin{equation}\label{2.6}
  \frac{d {\bf u}}{d t} + \mathcal{A} {\bf u} + \mathcal{B} p = {\bf f} ,
\end{equation}
\begin{equation}\label{2.7}
  \mathcal{B}^* {\bf u} = 0,
  \quad  0 < t \leq T,
\end{equation}
\begin{equation}\label{2.8}
  {\bf u}(0) = {\bf v} .
\end{equation}
For these operators in the space $\mathbf{L}_2(\Omega)$ we have 
\[
  \mathcal{A} = \mathcal{A}^* \geq \delta  \mathcal{E},
  \quad \delta = \delta(\Omega)  > 0 ,
\]
where $\mathcal{E}$ is the unit (identity) operator. 
Adjointness of operators $\mathcal{B} = \grad$ 
($\mathcal{B}: L_2(\Omega) \rightarrow \mathbf{L}_2(\Omega)$) and $\mathcal{B}^* = - \div$ 
($\mathcal{B}^*: \mathbf{L}_2(\Omega) \rightarrow L_2(\Omega)$)
follows from 
\[
  \int_{\Omega} {\bf u} \, \grad \, p \ d \, {\bf x} + 
  \int_{\Omega} \div \, {\bf u} \, p \ d \, {\bf x} = 0.
\]

For problem (\ref{2.6})--(\ref{2.8}) the following simple a priori estimate is valid
\begin{equation}\label{2.9}
  \| {\bf u} (t) \|^2 \leq \| {\bf v} \|^2 +
  \frac{1}{2 \delta} \int_{0}^{t} \| {\bf f} (\theta) \|^2 d \theta ,
\end{equation}
where $\| \cdot \|$ is the norm in $\mathbf{L}_2(\Omega)$.
Estimate (\ref{2.9}) will be for us a reference point when considering discrete problems. 

\section{Approximation in space}

In this work the main attention is paid to computational 
algorithms for the transition to a new time level, i.e. approximation in time. 
To construct discretization in time, operator-splitting schemes are used that allows to formulate a
problem for the pressure in the most natural way.
The problem of approximation in space is solved in the standard manner. 

There are employed various types
of grids: the non-staggered (collocated) grid where both the pressure and
velocity components are referred to the same points; next, partially staggered
(ALE-type) grid where the pressure is referred to an individual grid shifted in all space
directions on a one-half step from the basic grid where all velocity components
are defined; and finally, the staggered (MAC-type) grid where the pressure
is defined at the center points of grid cells whereas the velocity components
are referred to the corresponding faces of the cell. 

For simplicity we consider here uniform rectangular non-staggered grids. 
Problem (\ref{2.1})--(\ref{2.5}) is solved in a rectangle 
\[
  \Omega = \{ \ \mathbf{x} \ | \ \mathbf{x} = (x_1, x_2), 
  \ 0 < x_{\alpha} < l_{\alpha}, \ \alpha =1,2 \}.
\]
The approximate solution is calculated at the points of a uniform rectangular grid in $\Omega$:
\[
   \bar{\omega} = \{ \mathbf{x} \ | \ \mathbf{x} = (x_1, x_2),
   \quad x_\alpha = i_\alpha h_\alpha,
   \quad i_\alpha = 0,1,...,N_\alpha,
   \quad N_\alpha h_\alpha = l_\alpha\} 
\]
and let $\omega$ be the set of internal nodes 
($\bar{\omega} = \omega \cup \partial \omega$). 

For vector grid functions 
${\bf u}(\mathbf{x}) = 0, \ \mathbf{x} \in \partial \omega$
we define a Hilbert space ${\bf H} = {\bf L}_2({\omega})$ 
with the scalar product and norm 
\[
  ({\bf u},{\bf v}) = \sum_{{\bf x} \in \omega}
  {\bf u}({\bf x}) {\bf v}({\bf x}) h_1 h_2,
  \quad \|{\bf u}\| = ({\bf u},{\bf u})^{1/2} .
\]
The grid operator $A$ is taken in the form $A = - \nu \Delta_h$, 
where $\Delta_h$ is  the grid Laplace operator: 
\[
  \Delta_h y =
  - \frac{1}{h_1^2} (y(x_1+h_1,x_2) - 2y(x_1,x_2) + y(x_1-h_1,x_2))
\]
\[
  - \frac{1}{h_2^2} (y(x_1,x_2+h_2) - 2 y(x_1,x_2) + y(x_1,x_2-h_2)) .
\]
In ${\bf H}$ the operator $A$ is selfadjoint and positive definite: 
\begin{equation}\label{3.1}
  A = - \nu \Delta_h = A^* \geq \nu \delta_h E,
  \quad \delta_h = \sum_{\alpha = 1}^{2}
  \frac{4}{h^2_{\alpha}} \sin^2 \frac{\pi h_{\alpha}}{2 l_{\alpha}} .
\end{equation}
The pressure gradient is approximated by directed differences with an error $O(h)$.
We set $B = \grad_h$ at 
\begin{equation}\label{3.2}
  B p = \{ (B p)_1, (B p)_2 \},
  \quad \mathbf{x} \in \omega ,
\end{equation}
where  
\[
  (B p)_1 = \frac{1}{h_1} (p(x_1+h_1,x_2) - p(x_1,x_2)),
\]
\[
  (B p)_2 = \frac{1}{h_2} (p(x_1,x_2+h_2) - p(x_1,x_2)),
\]
The set of points for the pressure evaluation is denoted as
$\omega_p$ ($\omega_p \subset \bar{\omega}$).
For the grid divergence operator $B^* = -\div_h$ we have 
\begin{equation}\label{3.3}
  B^*{\bf u} = - \frac{1}{h_1} (u_1(x_1,x_2) - u_1(x_1-h_1,x_2))+
\end{equation}
\[
  - \frac{1}{h_2} (u_2(x_1,x_2) - u_2(x_1,x_2-h_2)) ,
  \quad \mathbf{x} \in \omega_p .  
\]
The adjointness  property of the grid gradient  and divergence operators 
is a consequence of the discrete equation 
\[
  \sum_{{\bf x} \in \omega}
  B p({\bf x}) {\bf u}({\bf x}) ({\bf x}) h_1 h_2 +
  \sum_{{\bf x} \in \omega_p}
  p({\bf x}) B^* {\bf u}({\bf x}) ({\bf x}) h_1 h_2 =0,
\]
which takes place on the set of vector grid functions 
${\bf u}(\mathbf{x}) = 0, \ \mathbf{x} \in \partial \omega$. 

In view of (\ref{3.1})--(\ref{3.3}) after the spatial approximation
of problem (\ref{2.6})--(\ref{2.8}) we obtain the following problem 
\begin{equation}\label{3.4}
  \frac{d {\bf u}}{d t} + A {\bf u} + B p = {\bf f} ,
\end{equation}
\begin{equation}\label{3.5}
  B^* {\bf u} = 0,
  \quad  0 < t \leq T,
\end{equation}
\begin{equation}\label{3.6}
  {\bf u}(0) = {\bf v} .
\end{equation}
For the solution of problem (\ref{3.4})--(\ref{3.6})
a priori estimate (\ref{2.9}) holds, where now $\| \cdot \|$ is
the norm in ${\bf H} = {\bf L}_2({\omega})$.

\section{Domain decomposition}

Let $\Omega $ be a combination of $p$ particular subdomains
$$
   \Omega =\Omega_1\cup \Omega_2\cup ...\cup \Omega_m \, .
$$
Particular subdomains can overlap one onto another. We shall construct the
schemes of decomposition where the solution at the new time level for the initial
problem is reduced to the sequential solution of problems in particular
subdomains.

Let us define functions for domain $\Omega $
\begin{equation}
\label{4.1}
   \eta_\alpha ({\bf x})=\left\{
       \begin{array}{c}
          >0, \quad {\bf x} \in \Omega_\alpha , \\
           0, \quad {\bf x} \notin \Omega_\alpha,
       \end{array}
   \quad \alpha =1,2, ..., m .\right.
\end{equation}
Generally, see for example, \cite{mathew2008domain,samarskii2002difference}, 
domain decomposition schemes for unsteady problems are based on the partition 
of unit for the region $\Omega$, where 
\[
   \sum_{\alpha =1}^{m} \eta_\alpha ({\bf x}) = 1,
   \quad {\bf x}\in \Omega \, .
\]
It is more convenient to use a somewhat different partition where 
\begin{equation}\label{4.2}
   \sum_{\alpha =1}^{m} \eta^2_\alpha ({\bf x}) = 1,
   \quad x\in \Omega \, .
\end{equation}

For the decomposition of  computational domain (\ref{4.1}), (\ref{4.2})
we consider the following additive representation of the identity operator 
$E$ in ${\bf H} = {\bf L}_2({\omega})$: 
\begin{equation}\label{4.3}
  E = \sum_{\alpha =1}^{m} \chi^2_\alpha,
  \quad \chi_\alpha = \eta_\alpha ({\bf x}) E,
  \quad {\bf x} \in \omega,
  \quad \alpha =1,2, ..., m .
\end{equation}
Taking into account (\ref{4.3}) we have 
\begin{equation}\label{4.4}
  {\bf u} = \sum_{\alpha =1}^{m} {\bf u}_\alpha,
  \quad {\bf u}_\alpha = \chi_\alpha {\bf u},
  \quad \alpha =1,2, ..., m .
\end{equation}

To formulate an appropriate system of equations for determining components 
of the solution ${\bf u}_\alpha, \ \alpha =1,2, ..., m$, 
we multiply both sides of equation (\ref{3.4}) by $\chi_\alpha$. 
This gives 
\begin{equation}\label{4.5}
  \frac{d {\bf u}_\alpha}{d t} + 
  \chi_\alpha A \sum_{\beta =1}^{m} \chi_\beta {\bf u}_\beta + B_\alpha p =
  {\bf f}_\alpha,
\end{equation}
where  
\[
  B_\alpha = \chi_\alpha B,
  \quad {\bf f}_\alpha = \chi_\alpha {\bf f},
  \quad \alpha =1,2, ..., m .
\]
Taking into account that 
\[
  B^*_\alpha = B^* \chi_\alpha ,
\]
equation (\ref{3.5}) in the new notation is written as 
\begin{equation}\label{4.6}
  \sum_{\alpha =1}^{m} B^*_\alpha {\bf u}_\alpha = 0,
  \quad  0 < t \leq T .
\end{equation}
The system of equations (\ref{4.5}), (\ref{4.6}) 
is supplemented by the initial conditions 
\begin{equation}\label{4.7}
  {\bf u}_\alpha(0) = {\bf v}_\alpha,
  \quad {\bf v}_\alpha = \chi_\alpha {\bf v},
  \quad \alpha =1,2, ..., m .
\end{equation}

For $U = \{{\bf u}_1, {\bf u}_2, ..., {\bf u}_m \}$, 
we define in the space ${\bf H}^m$ the norm and inner product as follows 
\[
  (U, V)_m =
  \sum_{\alpha =1}^{m} ({\bf u}_\alpha, {\bf v}_\alpha),
  \quad \|U\|_m = (U, U)_m^{1/2} .
\]
Let us multiply scalarly  individual equations 
(\ref{4.5}) by ${\bf u}_\alpha, \ \alpha =1,2, ..., m$  
and add them together. Next, multiply equation (\ref{4.6}) by $p$. 
Taking into account (\ref{4.4}), we obtain 
\[
  \frac{1}{2} \frac{d}{d t} \sum_{\alpha =1}^{m} 
  ({\bf u}_\alpha, {\bf u}_\alpha) +
  (A {\bf u},{\bf u}) =  \sum_{\alpha =1}^{m} 
  ({\bf f}_\alpha, {\bf u}_\alpha) .
\]
This implies the a priori estimate 
\begin{equation}\label{4.8}
  \|U\|_m^2 \leq \exp(t) \, \|V\|_m^2 +
  \int_{0}^{t} \exp(t-\theta) \, \|F(\theta)\|_m^2 d \theta 
\end{equation}
for problem (\ref{4.5})--(\ref{4.7}) with
$F = \{{\bf f}_1, {\bf f}_2, ..., {\bf f}_m \}$.

\section{Splitting scheme}

In the construction of domain decomposition schemes
we shall proceed from the scheme of splitting with respect to physical processes for the Cauchy 
(\ref{2.6})--(\ref{2.8}). 
We shall use a simple additive scheme componentwise splitting 
\cite{Marchuk:1990:SAD,0963.65091}.
Let ${\bf u}^n$ be the difference solution at the time moment 
$t^n=n\tau$,
where $\tau = T/N >0$ is the time-step. 
Let us separate out a particular stage connected with the pressure impact \cite{glowinski2003finite,guermond2006overview}.
Thus, in the first stage we have 
\begin{equation}\label{5.1}
  \frac{{\bf u}^{n+1/2} - {\bf u}^{n}}{\tau} + 
  A {\bf u}^{n+1/2} = {\bf f}^{n+1/2} .
\end{equation}
The pressure gradient is treated only in the second stage: 
\begin{equation}\label{5.2}
  \frac{{\bf u}^{n+1} - {\bf u}^{n+1/2}}{\tau} + B p^{n+1} = 0,
\end{equation}
\begin{equation}\label{5.3}
  B^* {\bf u}^{n+1} = 0 .
\end{equation}
Implementation of (\ref{5.2}), (\ref{5.3}) consists of two steps. 
In the first step we solve the following problem for the pressure
\[
  B^* B p^{n+1} = \frac{1}{\tau} \, B^* {\bf u}^{n+1/2},
\]
whereas in the second one we update the velocity: 
\[
  {\bf u}^{n+1} = {\bf u}^{n+1/2} - \tau B p^{n+1} .  
\]
Multiplying (\ref{5.1}) by ${\bf u}^{n+1/2}$, we obtain 
\[
  \|{\bf u}^{n+1/2}\|^2 \leq \|{\bf u}^{n}\|^2 +
  \frac{\tau}{\delta_h}  \|{\bf f}^{n+1/2}\|^2 .
\]
Similarly, from (\ref{5.2}) taking into account (\ref{5.3}) we have
\[
  \|{\bf u}^{n+1}\|^2 \leq \|{\bf u}^{n+1/2}\|^2 .
\]
Thus, we obtain the grid analog of estimate (\ref{2.9})
\[
  \|{\bf u}^{n+1}\|^2 \leq \|{\bf u}^{n}\|^2 +
  \frac{\tau}{2 \delta_h}  \|{\bf f}^{n+1/2}\|^2 .
\]
for difference scheme (\ref{5.1})--(\ref{5.3}).

For simplicity, we shall construct splitting schemes for problem 
(\ref{4.5})--(\ref{4.7}) by analogy with splitting scheme 
(\ref{5.1})--(\ref{5.3}).
The first half-step (viscous dissipation) is associated 
with the solution of equations 
\[
  \frac{d {\bf u}_\alpha}{d t} + 
  \chi_\alpha A \sum_{\beta =1}^{m} \chi_\beta {\bf u}_\beta =
  {\bf f}_\alpha,
  \quad \alpha =1,2, ..., m ,
  \quad t^n < t \leq t^{n+1/2}.
\]
Taking into account the fact that the numerical solution is implemented via solving
individual problems in the subdomains, 
the transition from time level $t^{n}$ to level $t^{n+1/2}$ 
can be realized as follows: 
\begin{equation}\label{5.4}
  \frac{{\bf u}_\alpha^{n+1/4} - {\bf u}_\alpha^{n}}{\tau} + 
  \chi_\alpha A \sum_{\beta =1}^{\alpha -1 }  \chi_\beta {\bf u}_\beta^{n+1/4} +
  \frac{1}{2}\chi_\alpha A \chi_\alpha {\bf u}_\alpha^{n+1/4} = {\bf f}_\alpha^{n+1/2},
\end{equation}
\[
  \alpha =1,2, ..., m ,
\]
\begin{equation}\label{5.5}
  \frac{{\bf u}_\alpha^{n+1/2} - {\bf u}_\alpha^{n+1/4}}{\tau} + 
  \frac{1}{2}\chi_\alpha A \chi_\alpha {\bf u}_\alpha^{n+1/2} +
  \chi_\alpha A \sum_{\beta = \alpha+1}^{m} \chi_\beta {\bf u}_\beta^{n+1/2} = 0,
\end{equation}
\[
  \alpha =1,2, ..., m .
\]
In view of (\ref{5.4}), (\ref{5.5}) in each subdomain 
$\Omega_\alpha, \ \alpha = 1,2, ..., m$ we must invert 
grid selfadjoint elliptic operator 
\[
  D_\alpha = E + \frac{1}{2} \chi_\alpha A \chi_\alpha
\]
for finding ${\bf u}_\alpha^{n+1/4}$ ($\alpha =1,2, ..., m$) and ${\bf u}_\alpha^{n+1/2}$
($\alpha =m,m-1, ..., 1$).
In this case, outside subdomains $\Omega_\alpha, \ \alpha = 1,2, ..., m$
there are used explicit calculations. 

Stability of scheme (\ref{5.4}), (\ref{5.5}) will be 
investigated in ${\bf H}^m$.
Consider the operator 
\begin{equation}\label{5.6}
  \mathbb{A} = \{ A_{\alpha \beta} \},
  \quad A_{\alpha \beta} = \chi_\alpha A \chi_\beta,
  \quad \alpha , \beta =1,2, ..., m .
\end{equation}
Taking into account (\ref{3.1}), (\ref{4.2}), we have 
$\mathbb{A} = \mathbb{A}^* \geq 0$ in ${\bf H}^m$.
Scheme \ref{5.4}), (\ref{5.5}) is based on 
using the triangular splitting 
\begin{equation}\label{5.7}
  \mathbb{A} = \mathbb{A}_1 + \mathbb{A}_2,
  \quad  \mathbb{A}_1 = \mathbb{A}^*_2 .
\end{equation}
Using notation (\ref{5.6}), (\ref{5.7}) we rewrite 
\ref{5.4}), (\ref{5.5})  in the form 
\begin{equation}\label{5.8}
  \frac{U^{n+1/4} - U^{n}}{\tau} + \mathbb{A}_1 U^{n+1/4} = F^{n+1/2} ,
\end{equation}
\begin{equation}\label{5.9}
  \frac{U^{n+1/2} - U^{n+1/4}}{\tau} + \mathbb{A}_2 U^{n+1/2} = 0 .
\end{equation}
Taking into account that $\mathbb{A}_\alpha \geq 0, \ \alpha =1,2$ 
in ${\bf H}^m$, 
for (\ref{5.9}) we immediately have 
\begin{equation}\label{5.10}
  \|U^{n+1/2}\|_m^2 \leq \|U^{n+1/4}\|_m^2 . 
\end{equation}
Multiplying (\ref{5.8}) by  $U^{n+1/4}$, we obtain 
\[
  \|U^{n+1/4}\|_m^2 \leq \|U^{n}\|_m^2 + 
  2 \tau (F^{n+1/2}, U^{n+1/4})_m .  
\]
For the last term on the right hand side we use the estimate 
\[
  2 \tau (F^{n+1/2}, U^{n+1/4})_m \leq 
  (1 - \exp(-\tau)) \, \|U^{n+1/4}\|_m^2 
  + \frac{\tau^2}{1 - \exp(-\tau))} \, \|F^{n+1/2}\|_m^2 .  
\]
This leads to the estimate 
\begin{equation}\label{5.11}
  \|U^{n+1/4}\|_m^2 \leq \exp(\tau) \, \|U^{n}\|_m^2 + 
  \tau \|F^{n+1/2}\|_m^2 .
\end{equation}

The second half-step results from the pressure and 
is connected with the system of equations 
\[
  \frac{d {\bf u}_\alpha}{d t} + B_\alpha p = 0,
  \quad \alpha =1,2, ..., m ,
\]
\[
  \sum_{\alpha =1}^{m} B^*_\alpha {\bf u}_\alpha = 0,
  \quad t^{n+1/2} < t \leq  t^{n+1}.
\]
Approximation in time for such systems 
were considered in \cite{vabishchevich2010additive}.
We shall use the additive scheme 
\begin{equation}\label{5.12}
  {\bf u}_\alpha^{n+1/2+\beta/2m} = {\bf u}_\alpha^{n+1/2+(\beta-1)/2m},
  \quad \beta \neq \alpha,
  \quad \beta =1,2, ..., m ,
\end{equation}
\begin{equation}\label{5.13}
  \frac{{\bf u}_\alpha^{n+1/2+\alpha/2m} - {\bf u}_\alpha^{n+1/2+(\alpha-1)/2m}}{\tau} + 
  B_\alpha p^{n+1/2+\alpha/2m} = 0,
\end{equation}
\begin{equation}\label{5.14}
  B_\alpha^* {\bf u}_\alpha^{n+1/2+\alpha/2m} = 0,
  \quad \alpha =1,2, ..., m .
\end{equation}
The implementation of additive scheme (\ref{5.12})--(\ref{5.14}) 
is conducted by analogy with scheme (\ref{5.2}), (\ref{5.3}).

For (\ref{5.12})--(\ref{5.14}) we have 
\[
  \|{\bf u}_\alpha^{n+1} \| \leq \|{\bf u}_\alpha^{n+1/2} \| ,
  \quad \alpha =1,2, ..., m 
\]
and therefore 
\begin{equation}\label{5.15}
  \|U^{n+1}\|_m^2 \leq \|U^{n+1/2}\|_m^2 . 
\end{equation}
Taking into account (\ref{5.10}), (\ref{5.11}) and (\ref{5.15})
we obtain the desired stability estimate of the additive 
operator-difference scheme (\ref{5.4}), (\ref{5.5}), (\ref{5.12})--(\ref{5.14})
\begin{equation}\label{5.16}
  \|U^{n+1}\|_m^2 \leq  \exp(\tau) \, \|U^{n}\|_m^2 + 
  \tau \|F^{n+1/2}\|_m^2 .
\end{equation}
Estimate (\ref{5.16}) is the grid analog of (\ref{4.8}) 
for the differential problem. 

This allows us to formulate the following main result. 

\begin{thm} 
The additive scheme of domain decomposition (\ref{5.4}), (\ref{5.5}), (\ref{5.12})--(\ref{5.14}) is unconditionally stable 
and estimate (\ref{5.16}) holds for the numerical solution.
\end{thm}

% \bibliographystyle{amsplain}
% \bibliography{DDM-NS}

\end{document}